\documentclass{birkart}
\usepackage{url,graphicx,color}
\usepackage{amssymb,amsthm,amsmath,paralist,umlaut}
\graphicspath{{EPS/}}

\theoremstyle{plain}
\newtheorem{lemma}{Lemma}[section]
\newtheorem{theorem}[lemma]{Theorem}

\newtheorem{proposition}[lemma]{Proposition}

\theoremstyle{definition}
\newtheorem{definition}[lemma]{Definition}

\theoremstyle{remark}

\newcommand\F{{\mathbb F}}
\newcommand\Q{{\mathbb Q}}
\newcommand\Z{{\mathbb Z}}
\newcommand\R{{\mathbb R}}

\newcommand\eps{\varepsilon}

\DeclareMathOperator{\aff}{aff}

\begin{document}
\title{Polyhedral surfaces of high genus} 
\author{Günter M. Ziegler}
\maketitle
\begin{abstract}
  The construction of the \emph{combinatorial} data for a surface with $n$
  vertices of maximal genus is a classical problem: The
  maximal genus $g=\lfloor\frac1{12}{(n-3)(n-4)}\rfloor$ was achieved
  in the famous ``Map Color Theorem'' by Ringel et al.\ (1968).
  We present the nicest one of Ringel's constructions,
  for the case $n\equiv7\bmod12$,
  but also an alternative construction, essentially due to Heffter (1898), 
  which easily and explicitly yields surfaces of genus $g\sim\frac1{16}n^2$.
  
  For \emph{geometric} (polyhedral) surfaces with $n$
  vertices the maximal genus is not known. 
  The current record is $g\sim n\log n$,
  due to McMullen, Schulz \& Wills (1983).
  We present these surfaces with a new construction: We
  find them in Schlegel diagrams of
  ``neighborly cubical $4$-polytopes,''
  as constructed by Joswig \& Ziegler (2000). 
\end{abstract}

\setcounter{section}{-1}%
\section{Introduction}

In the following we present constructions for surfaces that have
extremely and perhaps surprisingly high topological complexity (genus,
Euler characteristic) compared to their number of vertices.  We
believe that not only the resulting surfaces, but also the
constructions themselves are interesting and worth studying --- also
in the hope that they can be substantially improved).

\subsection{What is a surface?}

What do we mean by ``a surface''? This is not a stupid question, 
since combinatorialists, geometers, and topologists work
with quite different frameworks, definitions and concepts of surfaces,
and, as we will see, in the high-genus case it is not
clear that the various concepts coincide.

A \emph{topological} surface may be defined as a closed (compact,
without boundary), connected, 
orientable, Hausdorff, $2$-dimensional manifold. By adopting this
model, we already indicate that one could have worked in much greater
generality: Here we do not consider the non-orientable
case, we do not worry about manifolds with boundary, etc.

The \emph{combinatorial} version of a surface may be presented by listing the
\emph{faces} (vertices, edges, and, $2$-cells), 
and giving the necessary incidence information 
(for example, by specifying for each face the vertices in its boundary, 
in clockwise order according to the orientation).
Such combinatorial data must, of course, satisfy 
some consistency conditions if we are to be 
guaranteed that they do correspond to a surface.
Such conditions are easy to derive.

In the following, we will insist throughout that the 
combinatorial surface data we look at 
are regular (no identifications on the boundaries
of the cells), and they must satisfy the
\emph{intersection condition}: The intersection of 
any two faces is again a face (which may be empty).
This condition implies that any two vertices are connected by
at most one edge, and that any two $2$-faces have at most
two vertices in common (which must then be connected by an edge).%
\footnote{A combinatorial surface with the intersection condition
is called a ``polyhedral map'' in some of the 
discrete geometry literature; see Brehm \& Wills \cite{Breh}.}

\emph{Geometric} surfaces are embedded with flat faces in~$\R^3$.
Their faces are convex polygons, and we also require
that all these faces are simultaneously realized in~$\R^3$,
without intersections. 
Any such {geometric surface} yields a
combinatorial surface, which in turn yields a topological
manifold. 

\subsection{The \boldmath$f$-vector}

The $f$-vector of a combinatorial or geometric surface $S$ is the triple 
\[
f(S)\ := (f_0,f_1,f_2),
\]
where $f_0$ denotes the number of vertices,
$f_1$ is the number of edges, and $f_2$ is the number of
$2$-dimensional cells.

The $f$-vector contains a lot of information.  For example,
we can tell from the $f$-vector whether the surface is simplicial.
Indeed, one always has $3f_2\le2f_1$, by double-counting: Every face
has at least three edges, every edge lies in two faces.  Equality
$3f_2=2f_1$ holds if and only if every face is bounded by exactly
three edges, that is, for a triangulated (simplicial) surface.  

Similarly, we have $f_1\le\binom{f_0}2$, with
equality for a neighborly surface (with a complete graph),
which is necessarily simplicial.

\subsection{The genus}

The classification of the (orientable, closed, connected --- the
generality outline above) surfaces up to homeomorphism is well-known:
For each integer $g\geq 0$, there is exactly one topological type,
``the surface of genus $g$,'' which may be obtained by attaching  $g$
handles to the $2$-sphere $S^2$. 

The genus of a surface
may be defined, viewed, and computed in various different
ways, also depending on the model in which the surface is 
presented.

\emph{Topologically}, the genus may for example be obtained 
from a homology group, as $g=\frac{1}{2}\dim H_1(S_g; \Q)$. 
Alternatively, the 
genus may be expressed as the maximal number of disjoint,
non-separating, closed loops (this is the definition
given by Heffter \cite{heffter91:_ueber_probl_nachb}). 
It is also half the maximal number of non-separating loops that are
disjoint except for a common basepoint.

\emph{Combinatorially}, we
can compute the genus in terms of the Euler characteristic, 
$\chi(S_g)=2-2g=f_0-f_1+f_2$. So combinatorially the
genus is given by $g=1+\frac{1}{2}(f_1-f_0-f_2)\ge0$.

\subsection{The construction and realization problems}

Any combinatorial surface describes a topological space.
Conversely, any $2$-manifold can be triangulated,
but it is e.g.\ not at all clear how many vertices would
be needed for that. Thus we have the construction problem
for combinatorial surfaces:
\begin{description}
\item[Combinatorial construction problem]
    For which parameters $(f_0,f_1,f_2)$ are there
    combinatorial surfaces?
\end{description}
This is not an easy problem; in the triangulated case
of $2f_1=3f_2$ it is solved by Ringel's Map Color Theorem, discussed below.

Any geometric surface yields a combinatorial surface, but 
in the passage from combinatorial to geometric
surfaces, there are substantial open problems:
\begin{description}
\item[Geometric construction problem]
    For which parameters $(f_0,f_1,f_2)$ are there
    geometric surfaces?
\end{description}
This problem is much harder.
It may be factored into two steps, where the first
one asks for a classification or enumeration of
the combinatorial surfaces with the given parameters,
and the second one tries to solve the following
realization problem for all the combinatorial types:
\begin{description}
\item[Realization problem] 
  Which combinatorially given surfaces have geometric realizations?
\end{description}
 In general the answer to the Geometric construction problem does 
\emph{not} coincide with the answer for the Combinatorial
construction problem, that is,
the second step may fail even if the first one succeeds.
Let's look at some special cases:
\begin{compactitem}[~$\bullet$~]
  \item In the case of genus $0$, that is $f_1=f_0+f_2-2$,
   the construction problem was solved by Steinitz \cite{Stei3}: 
   The necessary and sufficient conditions both for
   combinatorial and for geometric surfaces 
   are $f_2\le2f_0-4$ and $f_0\le f_2-4$.
  \\
   By a second, much deeper, theorem by Steinitz \cite{Stei1}
   \cite{StRa} \cite[Lect.~4]{Z35}, every combinatorial surface of
   genus $0$ has a geometric realization in~$\R^3$, as the boundary of
   a convex polytope.
   This solves the realization problem for the case $g=0$.
  \item In the case of a simplicial torus, the
   possible $f$-vectors are easily seen to be $(n,3n,2n)$,
   for $n\ge7$.
A still pending, old conjecture of 
Grün\-baum \cite[Exercise 13.2.3, p.~253]{Gr1-2} states that every
triangulated torus (surface of genus~$0$, with $f$-vector
$(n,3n,2n)$) has a geometric realization in~$\R^3$.

  \item
On the other hand, there are combinatorial tori with
$f$-vector $(2n,3n,n)$, but \emph{none} of them 
has a geometric realization. Indeed, the condition 
$3f_0=2f_1$ means that the surface in question has a cubic
graph (all vertices of degree~$3$);
thus we are looking at the dual cell decompositions of
the simplicial tori. But none of them has a geometric
realization: Any geometric surface with a cubic graph
is necessarily convex --- that is, a $2$-sphere
(cf.\ \cite[Exercise 11.1.7, p.206]{Gr1-2}).

   \item  
Rather little is known about geometric surfaces of 
genus $g\ge2$:
Lutz \cite{lutz:_const} enumerated that there are
$865$ triangulated surfaces of genus $2$ on $10$
vertices, as enumerated by Altshuler. At least $827$ of
these have geometric realizations.
Specific examples of geometric surfaces of genus~$g\le4$ with
a minimal number of vertices were constructed by Brehm and
Bokowski \cite{BokowskiBrehm-genus3,BokowskiBrehm-genus4}.

  \item
There are also triangulated combinatorial surfaces that have
no geometric realization: 
Let's look at the $f$-vector $(12,66,44)$, which corresponds to
a neighborly surface of genus $6$ with $12$
vertices. Amos Altshuler has enumerated that there are exactly $59$
types of neighborly triangulations. One single one, number $54$, which
is particularly symmetric, was shown to be not geometrically realizable
by Bokowski \& Guedes de Oliveira \cite{BoGO2}. Thus,
$58$ possible types remain, and we do not know for any single one
whether it can be realized or not.
 \end{compactitem}
In general, it seems difficult to show 
for any given triangulated surface that
no geometric realization exists.
Besides the oriented matroid methods of Guedes de Oliveira \&
Bokowski, the obstruction theory set-up of Novik~\cite{Novik} 
and a linking-number approach of Timmreck~\cite{timmreck04}
have been developed in an attempt to do such non-realizability proofs.

In these lectures we look at families of combinatorial surfaces
whose genus grows quadratically in the number of vertices, such
as the neighborly triangulated surfaces on $n\gg7$ vertices,
where we \emph{think} that no geometric realizations exist,
but no general methods to prove such a general result seem to
be available yet.
And we present a construction for surfaces of genus $n\log n$,
which may be considered ``high genus'' in the category of
geometric surfaces, hoping that someone will be able to show
that this is good, or even best possible, or to improve upon it.

\section{Two combinatorial constructions}

Let us now look at a combinatorial surface with $f_0=n$ vertices. 
The following upper bound is quite elementary ---
the challenge is in the construction of examples that 
meet it, or at least get close.

\begin{lemma}
 A combinatorial surface with $n$ vertices has genus at most
\begin{equation}  \label{eq:g-bound}
  g\ \le\ \tfrac{1}{12}(n-3)(n-4).
\end{equation}
Equality can hold only for a triangulated surface that is neighborly,
which implies that $n$ is congruent to $0,3,4$ or $7\bmod12$.  
\end{lemma}

\begin{proof}
  Due to the intersection condition, any two vertices are connected by
  at most one edge, and thus $f_1\le\binom n2$.

In the case of a triangulated/simplicial surface, we have
$3f_2=2f_1$. With this, a simple calculation yields
\[
  g\ =  \ 1-\tfrac12(f_0-f_1+f_2) 
   \ =  \ 1-\tfrac12f_0+\tfrac16f_1
   \ \le\ 1-\tfrac12n+\tfrac16\tbinom n2
   \ =\ \tfrac1{12}(n-3)(n-4).
\]
 This holds with equality only if the surface is neighborly,
 and this can happen only if $\tfrac1{12}(n-3)(n-4)$ is an integer,
 that is, if $n\equiv 0,3,4$ or $7\bmod12$.  

If the surfaces is not simplicial, then it 
can be triangulated by introducing diagonals,
without new vertices, and without changing the genus.
However, this always results in triangulated surfaces with 
missing edges (diagonals that have not be chosen),
and thus in surfaces that do not achieve equality in (\ref{eq:g-bound}).
\end{proof}

The case of neighborly surfaces is indeed very interesting, and
has received a lot of attention. In particular, it occurred first in
connection with (a generalization of) the four color problem: Its
analog on surfaces of genus $g>0$, known as the ``Problem der
Nachbargebiete,'' the problem of neighboring countries, is solved by
exhibiting of a maximal configuration of ``countries'' that are
pair-wise adjacent. If one draws the dual graph to such a
configuration, then this will yield a triangulation of the surface
(Kempe 1879 \cite{kempe79}; Heffter 1891 \cite{heffter91:_ueber_probl_nachb}). 
As the ``thread
problem'' (Fadenproblem) the question was presented in the famous book
by Hilbert \& Cohn-Vossen \cite{HilbertCohnVossen}. 

The case $n=4$ is trivial (realized by the tetrahedron); the first
interesting case is $n=7$, where a combinatorially-unique configuration exists,
the simplicial ``Möbius Torus'' on $7$ vertices
\cite{moebius86:_mitth_nachl}. We will look at it below.
Möbius' triangulation was rediscovered a number of times, realized by
Császár, and finally exhibited in the Schlegel diagram of a
cyclic $4$-polytope on $7$ vertices,
by Duke \cite{duke70:_geomet} and Altshuler \cite{Alts}. 
For the other neighborly cases, $n\ge12$, no realizations
are known.

When $n$ is not congruent to $0, 3,4,$ or $7$ the maximal genus of a
surface on $n$ vertices if of course smaller than the bound given
above, but it could be just the bound rounded down, and indeed it is.

\begin{theorem}\label{thm:MapColorTheorem}
  [Ringel et al.\ (1968); see \cite{ringel74:_map_color_theor}]
  For each $n\geq 4$, $n\neq9$,
  there is a (combinatorial) $n$-vertex surface of genus
\[\textstyle
g_{\max}\ =\ \left\lfloor{\frac{(n-3)(n-4)}{12}}\right\rfloor.
\]
  
\end{theorem} 

In his 1891 paper, Heffter \cite[\S3]{heffter91:_ueber_probl_nachb}
proved this theorem for $n\le 12$; in particular, in doing this
he introduced some of the basic concepts and notation, and thus ``set the
stage.'' 
From then, it needed another 77 years to complete the proof
of Theorem~\ref{thm:MapColorTheorem}.
The full proof is complicated, with intricate combinatorial arguments
divided into twelve cases (according to $n\mod  12$) and a
number of ad-hoc constructions needed for sporadic cases of ``small $n$.'' In
the following we will sketch Ringel's construction for the nicest of
the twelve cases, the case of $n\equiv 7\bmod  12$.  This is the
only case where we can get a surface with a cyclic symmetry, according
to Heffter, and in fact we do!  (This special case was
first solved by Ringel in 1961, but our exposition follows his book
from 1973, to which we also refer for the other eleven cases.)
Then we also present a second construction, based on a paper by
Heffter from 1898 \cite{heffter98:_ueber_grupp_nachb}: 
This produces surfaces that are not quite
neighborly, but they still do have genus that grows quadratically with
the number of vertices. Moreover, this construction is very conceptual
and explicit. For simplicity we will give a simple combinatorial
description, but indeed one may note that the surface has a 
$\Z_q$-action whose quotient is the ``perfect''
cellulation with just one vertex and one $2$-cell, and thus the
surface we get arises as an abelian covering from the perfect
cellulation, where opposite edges of a $4g$-gon are identified.

\subsection{A neighborly triangulation for \boldmath$n\equiv7\bmod  12$}

It was observed already by Heffter that 
a combinatorial surface is completely determined if we
label the vertices, and for each vertex describe the cycle
of its neighbors (in counter-clockwise/orientation order).

Thus, for example, a ``square pyramid'' (a $2$-sphere with
$5$ vertices, consisting of 
one quadrilateral and four triangles, is given by a 
\emph{rotation scheme} of the form
\begin{eqnarray*}
  0&:& (1,2,3,4)\\[-4pt]
  1&:& (0,4,2)\\[-4pt]
  2&:& (0,1,3)\\[-4pt]
  3&:& (0,2,4)\\[-4pt]
  4&:& (0,3,1)\\[-20mm]&&\hspace{10mm}\qquad
\input{EPS/delta1.pstex_t}
\end{eqnarray*}
which says that $1,2,3,4$ are the neighboring vertices, in cyclic
order, for vertex~$0$, etc.
In particular, we could have written $(2,3,4,1)$ instead of
$(1,2,3,4)$, since this denotes the same cyclic permutation.
Some checking is needed, of course, to see whether 
some scheme of this form actually describes a surface 
that satisfies the intersection condition.

In the case of a triangulated surface, the corresponding 
consistency conditions are rather easy to describe.
Indeed, if $j,k$ appear adjacent in the cyclic list
of neighbors to a vertex $i$, then this means that
$[i,j,k]$ is an oriented triangle of the surface ---
and thus $k,i$ have to be adjacent in this order in
the cycle of neighbors for $j$, and similarly $i,j$ have to appear
in the list for~$k$.

\begin{figure}[ht]
\begin{center}
\input{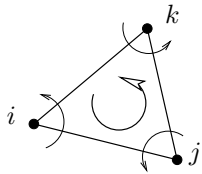}
\end{center}
\caption{Reading off data of the rotation scheme from a triangle in an
oriented surface.}
\label{figure:deltastar}
\end{figure}

Thus in terms of the rotation scheme, 
the \emph{triangulation condition} (which Ringel calls the 
``rule $\Delta^*$'') says that if the row for vertex $i$ reads
\begin{eqnarray*}
i &:& (\ \ldots\ \ldots\ \ldots\ ,j,k,\ \ldots\ \ldots\ \ldots\ )
\end{eqnarray*}
then in the rows for $j$ and $k$ we have to get
\begin{eqnarray*}
j &:& (\ \ldots\ \ldots\ ,k,i,\ \ldots\ \ldots\ \ldots\ \ldots\ )\\[-4pt]
k &:& (\ \ldots\ \ldots\ \ldots\ \ldots\ ,i,j,\ \ldots\ \ldots\ ).
\end{eqnarray*}

We want to construct triangulated surfaces
with a cyclic automorphism group $\Z_n$ ---
so the scheme for one vertex would yield all others by addition
modulo~$n$.
Unfortunately, this is possible \emph{only} for 
$n=4,5,6$ and for $n\equiv7\bmod12$, according to 
Heffter \cite[\S4]{heffter91:_ueber_probl_nachb}.

For example, for $n=7$ there is such a surface, the
Möbius torus \cite{moebius86:_mitth_nachl}, given by
\begin{eqnarray*}
0 &:& (1,3,2,6,4,5)\\[-4pt]
1 &:& (2,4,3,0,5,6)\\[-4pt]
2 &:& (3,5,4,1,6,0)\\[-4pt]
3 &:& (4,6,5,2,0,1)\\[-4pt]
4 &:& (5,0,6,3,1,2)\\[-4pt]
5 &:& (6,1,0,4,2,3)\\[-4pt]
6 &:& (0,2,1,5,3,4).\\[-30mm]&&\hspace{30mm}\qquad
\input{EPS/moebius2.pstex_t}
\end{eqnarray*}
Here the first row determines all others by addition modulo~$n$.

Now let's assume we have a rotation scheme for
a triangulated surface with $\Z_n$ automorphism group. 
If the row for vertex $0$ reads
\begin{eqnarray*}
0 &:& (\ \ldots\ \ldots\ \ldots\ ,j,k,\ \ldots\ \ldots\ \ldots\ )
\end{eqnarray*}
then the triangulation condition, rule $\Delta^*$, yields that
\begin{eqnarray*}
j &:& (\ \ldots\ \ldots\ ,k,0,\ \ldots\ \ldots\ \ldots\ \ldots\ )\\
k &:& (\ \ldots\ \ldots\ \ldots\ \ldots\ ,0,j,\ \ldots\ \ldots\ )
\end{eqnarray*}
and then the $\Z_n$-automorphism implies (subtracting $j$ resp.~$k$) that
\begin{eqnarray*}
0 &:& (\ \ldots\ \ldots\ ,k{-}j,-j,\ \ldots\ \ldots\ \ldots\ )\\
0 &:& (\ \ldots\ \ldots\ \ldots\ ,-k,j{-}k,\ \ldots\ \ldots\ ).
\end{eqnarray*}
In other words, if in the neighborhood of~$0$, we have that
``$k$ follows $j$,'' then also 
``$-j$ follows $k{-}j$,'' and 
``$j{-}k$ follows $-k$''
(where all vertex labels are interpreted in $\Z_n$, that is, 
modulo~$n$).

\begin{figure}[ht]
\begin{center}
\input{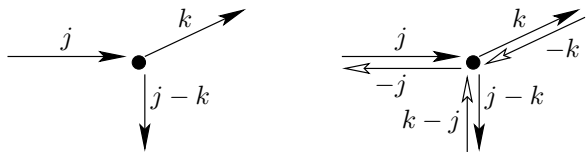}
\end{center}
\caption{The flow condition, ``Kirchhoff's law.'' 
The left figure shows how a flow of size $j$ is split into
two parts. In the right figure the reversed arcs have been added:
This contains the same amount of information, but the
data can be read off more directly.}
\label{figure:kirchhoff}
\end{figure}

The condition that we have thus obtained
can be viewed as a flow condition (a ``Kirchhoff law'') in
a cubic graph:
The cyclic order in the neighborhood of~$0$ can be 
derived from a walk in an edge-labelled graph,
whose edge labels satisfy a flow condition --- 
see Figure~\ref{figure:kirchhoff}.

\begin{figure}[ht]
\begin{center}
\input{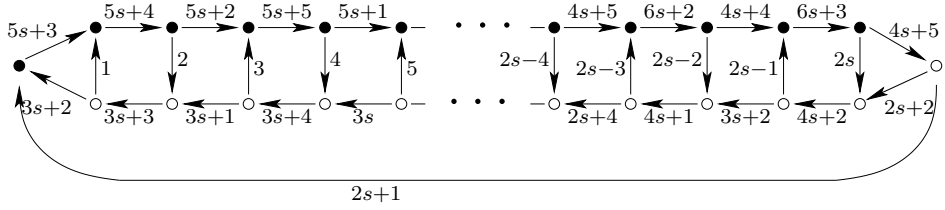}
\end{center}
\caption{A network for a neighborly surface with $n=12s+7$ vertices.} 
\label{figure:kirchhoff2}
\end{figure}

Thus in order to obtain a valid ``row~$0$'' we have to produce
a cyclic permutation of $1,2,\dots,n-1$ that can be read off
from a flow (circulation) in a cubic graph.
Ringel's solution for this in the case $n\equiv 7\bmod12$ is given by
Figure~\ref{figure:kirchhoff2}.

It is based on writing $\Z_n=\{0,\pm1,\pm2,\dots,\pm(6s+3)\}$.
The figure encodes the full construction:
It describes a cubic graph with $2s+4$ vertices and $6s+3$ edges,
where
\begin{compactitem}[--]
  \item each edge label from $\{1,2,\dots,6s+3\}$ occurs exactly once,
  \item at each vertex, the flow condition is satisfied (modulo $n$).
  \end{compactitem}
Now the construction rule is the following: Travel on this graph, 
\begin{compactitem}[--]
  \item at each 
    black vertex $\bullet$ turn ``left,''
    to the next arc in clockwise direction, at each
    white vertex $\circ$ turn ``right,''
    to the next arc in counterclockwise direction, and
  \item record the label of each edge traversed in arrow direction,
   resp.\ the negative of the label if traversed against arrow
   direction,
\end{compactitem}
The main claim to be checked is that this prescription
leads to a single cycle in which each edge is traversed
in each direction exactly once, so each value in
$\{\pm1,\pm2,\dots,\pm(6s+3)\}$ occurs exactly once.
For example, if we start at the arrow labelled ``$1$,'' then the
sequence we follow will start
\[
1,\ -(5s+3),\ -(3s+2),\ -(3s+3),\ -(3s+1),\ -(3s+4),\ -3s,\ -(3s+5),\ 
\ldots
\]
This is the first line of the rotation diagram for 
Ringel's neighborly surface with $n={12s+7}$ vertices. 

Note: \emph{any} cyclic order yields a surface, 
but we need to control the intersection property,
and the genus, e.g.\ by enforcing the triangle condition.
On the other side, if we just take a random permutation (cyclic
order), then this yields a very interesting model of a random
surface of random genus.
See Pippenger \& Schleich \cite{pippenger:_topol_char}
for a current discussion of such models.

\subsection{Heffter's surface and a triangulation}

Here is a much simpler construction, which
yields a not-quite neighborly surface.
The underlying remarkable cellular surface was first
discovered by Heffter \cite{heffter98:_ueber_grupp_nachb}, 
much later rediscovered by Eppstein et al.\ \cite{Z80}.
(See also Pfeifle \& Ziegler \cite{Z85}.)

Let $q=4g+1$ be a prime power with $g\ge1$ (one can find
suitable primes $q$, or simply take $q=5^r$).
The one algebraic fact we need is that there is
a finite field $\F_q$ with $q$ elements, and that
the multiplicative group $\F_q^*=\F_q\setminus\{0\}$ is cyclic 
(of order $q-1=4g$),
that is, there is a generator $\alpha\in\F_q^*$ such that
$\F_q^*=\{\alpha,\alpha^2,\alpha^3,\dots,\alpha^{q-1}\}$,
with $\alpha^{q-1}=\alpha^{4g}=1$.
In particular, we get $\alpha^{2g}=-1$.
For example, for $g=3$ and $q=13$ we may take
$\alpha=2$, with
$(1,\alpha,\alpha^2,\dots\ )=(1, 2, 4, 8, 3, 6, 12, 11, 9, 5, 10, 7)$.

\begin{figure}[ht]
\begin{center}
\input{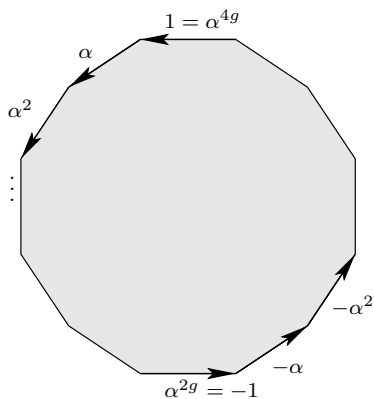}
\end{center}
\caption{Identifying the opposite edges of a  $4g$-gon
we obtain a perfect cellulation of~$S_g$: All vertices are 
identified.}
\label{figure:heffter1}
\end{figure}

For any $g\ge1$, a \emph{perfect} cellulation of~$S_g$
is obtained from a $4g$-gon by identifying opposite edges
in parallel. 
In the prime power case, a combinatorial description for this
is as follows: Label the directed edges of the $n$-gon by
$1,\alpha,\alpha^2,\alpha^3,\dots$ in cyclic order,
and identify the antiparallel edges labelled $s$ and $-s$.
(Compare Figure~\ref{figure:heffter1}.)

The resulting cell complex has the 
$f$-vector $f=(1,2g,1)$. It is \emph{perfect} in the sense of
Morse theory since this is also
the sequence of Betti numbers (ranks of the homology groups).
However, this cell decomposition is not ``regular'' in the
sense that there are identifications on the boundaries of cells:
All the ends of the edges are identified, and there are lots of
identifications on the boundary of the $2$-cell.

\begin{figure}[ht]
\begin{center}
\input{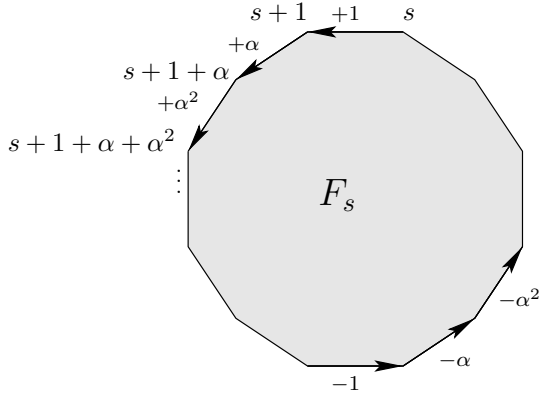}
\end{center}
\caption{One of the $q$ $2$-cells, and the labelling of
its $q-1=4g$ vertices (by elements of~$\F_q$).}
\label{figure:heffter2}
\end{figure}

Now we explicitly write down a $q$-fold ``abelian covering'' 
of this perfect cellulation. It
has both its vertices and its $2$-cells indexed by $\F_q$:
The surface has $q$-cell faces $F_s$, for $s\in\F_q$.
The vertices will also be labelled by the $q$ elements of~$F_q$.
Namely, the face $F_s$ should have vertices
\[
s,\ 
s+1,\ 
s+1+\alpha,\ 
s+1+\alpha+\alpha^2,\ \ \ldots\ \ ,
s+1+\alpha+\dots+\alpha^{4g-1},
\]
in cyclic order (as indicated by Figure~\ref{figure:heffter2}),
that is,
\[
s+\sum_{i=0}^{k-1} \alpha^i\ \ =\ \ 
s+\frac{\alpha^k-1}{\alpha-1}\qquad
\textrm{for }0\le k<4g-1.
\]
For each face $F_s$ this yields $q-1$ distinct values/vertices: 
$\alpha^k$ takes on every value except for $0$, and thus
$s+\frac{\alpha^k-1}{\alpha-1}$
yields all elements of~$\F_q$ except for $s+\frac{-1}{\alpha-1}$.
(Explicit worked out examples, for $q=5$ resp.\ for $9=9$,
can be found in \cite{Z80} and~\cite{Z85}.)

Now we have to verify that this prescription does indeed give
a surface: For this, check that each vertex comes to lie
in a cyclic family of $q-1$ faces.

We thus have a quite remarkable combinatorial structure:
The cellular surface $\widetilde S_g$ has $q$ vertices and $q$ faces;
the vertices have degree $q-1$, the faces have $q-1$ neighbors.
Thus the graph of the surface is complete (each vertex adjacent
to every other vertex), and so is the dual graph
(each face adjacent to every other one).
Moreover, the surface is self-dual, that is, isomorphic to the
dual cell decomposition.

With all the combinatorial facts just mentioned, we have
in particular computed the $f$-vector of the surface: It is
\[
f(\widetilde S_g)\ =\ (q,2gq,q)\ =\ (q,\tbinom q2,q).
\]
Thus we have an orientable surface with Euler characteristic 
$2q-\binom q2=2-2g$, and genus $g=\frac12\binom q2-q+1$.

Moreover, ``by construction'' the surface is very symmetric: First,
there clearly is an $\F_q$-action by addition; and if we mod
out by this action, we recover the original, ``perfect'' cell
decomposition $\widetilde S_g/\Z_q\cong S_g$ with one face.
But also the multiplicative group $\F_q^*$ acts by multiplication, with
$\alpha\cdot(s+   1   + \dots+\alpha^{i-1} )= 
      \alpha s+\alpha + \dots+\alpha^i      = 
  (\alpha s-1)+   1 +\alpha+\dots+\alpha^i$.
Thus the action is given by
\[F_s\mapsto F_{\alpha s-1}, \qquad
v_s^i\mapsto v_{\alpha s-1}^{i+1}.
\]
The full symmetry group of~$S_g$ is a 
``metacyclic group'' with $q(q-1)$ elements.

The surface $\widetilde S_g$ is a regular cell complex,
but it does not satisfy the intersection condition:
 Any two $2$-cells intersect in $q-2$ vertices 
(since each $2$-cell has $q-1$ vertices, that is, all 
of them except for one).

Thus we triangulate $\widetilde S_g$, by stellar subdivision of the $2$-cells:
Then we have $n=2q$ vertices, $q$ of degree $q-1$, and $q$ of degree $2q-2$.
Furthermore, there are $f_1=3\binom n2$ edges,
namely $\binom q2$ ``old'' ones and $q(q-1)$ ``new ones''
introduced by the $q$ stellar subdivisions. Furthermore, we
have now $q(q-1)$ triangle faces, which yields an $f$-vector
\[
f\ =\ \big(\, 2q,\, 3\tbinom q2,\, 2\tbinom q2\,\big),
\]
and hence 
\[
g\ =\ 1+\tfrac12(2q-3\tbinom q2+2\tbinom q2)
 \ =\ \frac{n^2-10n+16}{16}.
\]
So for these simplicial surfaces, for which we have a
completely explicit and very simple combinatorial description,
the genus is quadratically large in the number
of vertices, $g\sim \frac{n^2}{16}$,
but they don't quite reach the value
$g\sim\frac{n^2}{12}$ of  neighborly surfaces.

\subsection*{Conclusion}

So what is the moral? The moral is that using combinatorial
constructions, we do obtain triangulated surfaces whose genus grows
quadratically with the number of vertices. To find the constructions
for surfaces with
the exact maximal genus is very tricky, and certainly one would hope
for simpler and more conceptual descriptions/constructions,
but combinatorial surfaces whose genus grows quadratically with
the number of vertices are quite easy to get.

\section{A geometric construction}

Any smooth surface embedded or immersed in
$\R^3$, equipped with a generic ``height'' function,
as studied by Morse theory, conforms to a chain of
inequalities
\[
g=\dim H_1(S)<\dim H_{*}(S)\le \#\textrm{ critical points}.
\]
If we think of a simplicial/polyhedral surface in $\R^3$ as an
approximation to a smooth surface, then we would also use a linear
objective function (as a Morse function), might conclude that all the
critical points should certainly be at the vertices, and thus the
genus $g$ cannot be larger than the number of vertices for an embedded
(or immersed) surface.

\emph{However}, in the case of high genus the approximation of a smooth
surface by a simplicial surface is not good, it is very coarse, and
the critical points induced by a linear function on a simplicial
surface certainly do not satisfy the Morse condition of looking like
quadratic surfaces. And indeed, the result suggested by our argument
is far from being true. It was disproved by McMullen, Schulz \& Wills
\cite{mcmullen83:_polyh_e},
who in 1983 constructed ``polyhedral $2$-manifolds in
$E^3$ with unusually high genus'': They produced sequences both of
simplicial and of quad-surfaces, whose genus grows like $n\log n$ in
the number of vertices.

In the following we will give a simple combinatorial
description of ``their'' quad-surfaces $Q_m$, and describe an
explicit, new geometric construction for them,
which is non-inductive, yields explicit coordinates, and ``for free''
even yields a cubification of the convex hull of the surface without
new vertices. We obtain this by putting together (simplified versions
of) several recent constructions: Based on intuition from Amenta \&
Ziegler \cite{Z51a}, a simplified construction of the neighborly
cubical polytopes of Joswig \& Ziegler \cite{Z62}, which are connected
to the construction of high genus surfaces via Babson, Billera \&
Chan \cite{BBC} and observations of Schröder \cite{schroeder04}.  The
constructions as presented here can be generalized and extended quite
a bit, which constitutes both recent work as well as promising and
exciting directions for further research. See e.g.\ Ziegler
\cite{Z97}.

The construction in the following will be in five parts: 
\begin{compactenum}[1.]
\item Combinatorial description of the surface as the mirror-surface of
  the $n$-gon, embedded into the $n$-dimensional standard cube,
\item construction of a deformed $n$-cube,
\item definition and characterization of faces that are ``strictly
  preserved'' under a polytope projection,
\item identification of some faces of the deformed $n$-cube above that
  are strictly preserved under projection to $\R^4$, and
\item putting it all together, and obtaining the desired surfaces via
  Schlegel diagrams.
\end{compactenum}

\subsection{Combinatorial description}\label{subsec:combinatorialdesc}

The surface $Q_m$ is most easily described as a 
subcomplex of the $m$-dimensional cube~$C_m=[0,1]^m$.

Any nonempty face of~$C_m$ consists of those points in~$C_m$
for which some coordinates are fixed to be~$0$,
others are fixed to be~$1$, and the rest are left free to
vary in~$[0,1]$.
Thus there is a bijection of the non-empty faces 
with $\{0,1,{*}\}^m$.

\begin{definition}
  For $m\ge3$,
  the quad-surface $Q_m$ is given by all the faces of~$C_m$
  for which only two, cyclically-successive coordinates may be
  left free.
\end{definition}

Thus the subset $|Q_m|\subset[0,1]^m$ consists of all points
that have at most two fractional coordinates --- and if there
are two, they have to be either adjacent, or they have to be
the first and the last coordinate.
(This description perhaps first appeared in Ringel
\cite{ringel55:_ueber_probl_wuerf_wuerf}.)
In particular, $Q_3$ is just the boundary of the unit $3$-cube.

Let's list the faces of~$Q_m$: 
These are \emph{all} the $f_0(Q_m)=2^m$ vertices of the 
$0/1$-cube, encoded by $\{0,1\}^m$;
then $Q_m$ contains \emph{all} the $f_1(Q_m)=m2^{m-1}$ 
edges of the $m$-cube, corresponding to strings with exactly
one $*$ and $0/1$-entries otherwise.
And finally we have $f_2(Q_m)=m2^{m-2}$ quad faces,
corresponding to strings with two cyclically-adjacent $*$s
and $0/1$s otherwise.

\begin{figure}[ht]
\begin{center}
\input{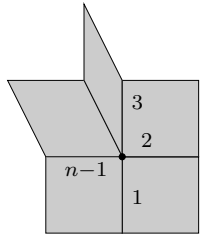}
\end{center}
\caption{The star of a vertex in $Q_m$}
\label{figure:Q_m-vertexstar}
\end{figure}

Why is this a surface? This is since all the
vertex links are circles. Indeed, if we look at any vertex, then we see  
in its star the $m$ edges emanating from the vertex, and the
$n$ square faces between them, which connects them in the
cyclic order, as in Figure~\ref{figure:Q_m-vertexstar}.

\begin{figure}[ht]
\begin{center}
\input{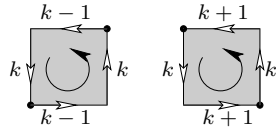}
\end{center}
\caption{The orientation of~$Q_m$ described in the text; here the
  black vertices are the ones with an even sum of coordinates.}
\label{figure:Q_m-orientable}
\end{figure}

It is similarly easy to see that the surface is indeed
orientable: An explicit orientation is obtained by dictating
that in the boundary of any $2$-face for which the fractional coordinates
are $k-1$ and $k$ (modulo~$m$), the edges with a fractional
$(k-1)$-coordinate should be oriented from the even-sum vertex to
the odd-sum vertex, while the edges corresponding to a fractional $k$-th
coordinate are oriented from the odd-sum vertex to the even-sum
vertex (cf.~Figure~\ref{figure:Q_m-orientable}).

Thus, $Q_m$ is an orientable polyhedral surface,
realized geometrically in $\R^m$ as a subcomplex of~$C_m$.
Its Euler characteristic is
\[
\chi(Q_m)\ =\ 2^m-m2^{m-1}+m2^{m-2}\ =\ (4-2m+2)2^{m-2}
\]
and thus with $g=1+\frac12\chi$ and $n:=f_0(Q_m)=2^m$ the genus is
\[
g(Q_m)\ =\ 1+(m-4)2^{m-3}
\ =\ 1+\tfrac n8\log\tfrac n{16}\ =\ \Theta(n\log n).
\]
So we are dealing with a $2$-sphere for $m=3$,
with a torus for $m=4$, while for $m=5$ we already
get a surface of genus $5$.
There are also simple recursive ways to describe the surface $Q_m$,
as given by McMullen, Schulz \& Wills in their original paper
\cite{mcmullen83:_polyh_e}.
The combinatorial description here is a special case of the
``mirror complex'' construction of Babson, Billera \& Chan~\cite{BBC},
which from any simplicial $d$-complex on $n$ vertices
produces a cubical $(d+1)$-dimensional subcomplex of the
$n$-cube on $2^n$ vertices and the given complex in all the vertex links:
Here we are dealing with the case of~$d=1$, where the simplicial
complex is a cycle on $n$ vertices.

\begin{figure}[ht]
\begin{center}
\includegraphics[height=77mm]{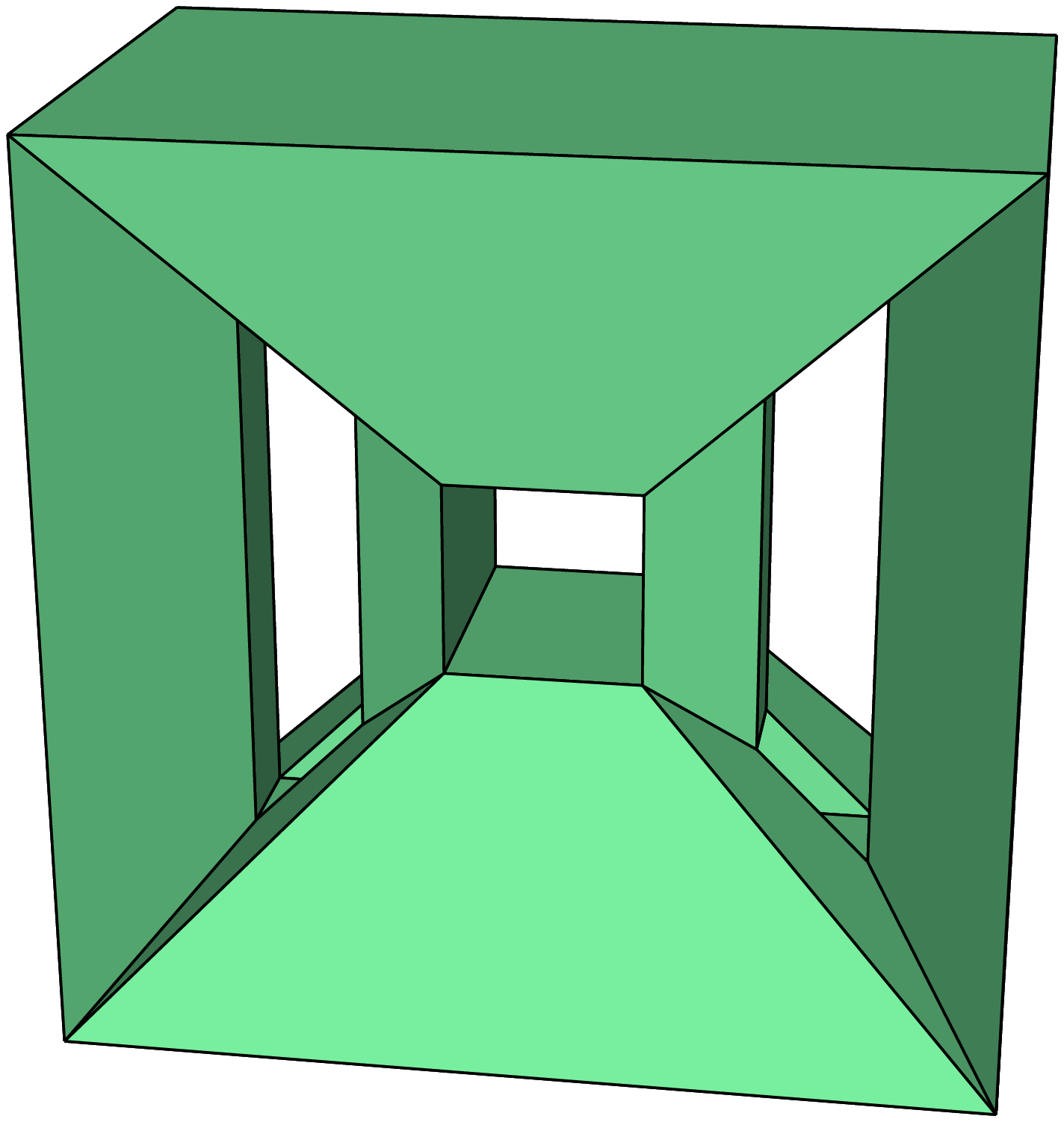}
\end{center}
\vspace{-4mm}
\caption{The surface $Q_5$ of genus~$5$, realized in~$\R^3$
(\texttt{polymake}/\allowbreak \texttt{javaview} graphics by Thilo Schröder)}
\label{figure:genus5surface}
\end{figure}

\subsection{Construction of a deformed \boldmath$m$-cube}
\label{subsec:constructdeformedcube}

In the last section, we have described a surface $Q_m$
as a subcomplex of the standard orthogonal $n$-cube $C_m=[0,1]^m$,
and thus as a polyhedral complex in~$\R^m$.
If we take any other realization of the $m$-cube, 
then this yields a corresponding realization of our surface as a subcomplex.
The object of this section is to describe an entirely explicit
``deformed'' cube realization $D_m^\eps$, which contains the
surface $Q_m^\eps$  as a subcomplex. Here it is.

\begin{definition}
  For $m\ge4$ and $\eps>0$ let
  $D_m^\eps$ be the set of all points $x\in\R^m$ that satisfy the 
  linear system of $2m$ linear inequalities
  \begin{equation}
    \label{eq:defcube}
\left(\begin{array}{r@{\hspace{5pt}}r@{\hspace{5pt}}r@{\hspace{5pt}}r@{\hspace{5pt}}r@{\hspace{5pt}}r@{\hspace{2pt}}|@{\hspace{2pt}}r@{\hspace{5pt}}r@{\hspace{5pt}}r@{\hspace{5pt}}r}
\pm\eps\\
  2 &\pm\eps\\
 -7 &  2 &\pm\eps\\
  7 & -7 &  2 &\pm\eps\\
 -2 &  7 & -7 &  2 &\pm\eps\\
    & -2 &  7 & -7 &\cdot&\cdot \\
    &    & -2 & \cdot  & -7 &  2 &\pm\eps\\
    &    &    & \cdot  &  7 & -7 &  2  &\pm\eps\\
    &    &    &        & -2 &  7 & -7 &  2 &\pm\eps\\
    &    &    &        &    & -2 &  7 & -7 & 2 &\pm\eps
  \end{array}\right).
\left(\begin{matrix}
x_1\\x_2\\x_3\\x_4\\\cdot\\\cdot\\\cdot\\x_{m-2}\\x_{m-1}\\x_m
  \end{matrix}\right)
\le
\left(\begin{matrix}
b_1\\b_2\\b_3\\b_4\\\cdot\\\cdot\\\cdot\\b_{m-2}\\b_{m-1}\\b_m
  \end{matrix}\right)
  \end{equation}
\end{definition}

This defines a polytope with the combinatorics of~$C_m$,
if $\eps>0$ is small enough, and if the
sequence of right-hand side entries $b_1,b_2,b_3,\dots$ grows fast enough.
The following lemma gives concrete values ``that work.''

\begin{lemma}
  For $0<\eps<\frac12$ and $b_k=(\frac6\eps)^{k-1}$,
  the set $D_m^\eps$ is combinatorially equivalent to the 
  $n$-cube.
\end{lemma}

\begin{proof}
The $k$-th pair of inequalities from (\ref{eq:defcube})
may be written as
\begin{equation}
  \label{eq:kth-ineq}
  |x_k|\ \le\ \frac1\eps
  (b_k - 2x_{k-1} + 7x_{k-2} - 7x_{k-3} +2x_{k-4}),
\end{equation}
with $x_0\equiv x_{-1}\equiv x_{-2}\equiv x_{-3}:\equiv 0$.
So if the $x_{k-1},x_{k-2},\dots$ are bounded, and
$b_k$ is guaranteed to be larger than 
\[  L_k \ :=\ 2|x_{k-1}|+7|x_{k-2}|+7|x_{k-3}|+2|x_{k-4}|,\]
then the right-hand side of (\ref{eq:kth-ineq}) is 
strictly positive, and the $x_k$ is bounded again.
In this situation, we find that the two inequalities represented
by (\ref{eq:kth-ineq}) cannot be simultaneously satisfied with 
equality, but any one of them can.
Thus, inductively we get that the first $2k$ inequalities of the
system (\ref{eq:defcube}) define a $k$-cube (in the first $k$
variables).

With the concrete values as suggested by the lemma, we
verify by induction that $|x_k|\le\frac13(\frac6\eps)^k$.
Indeed, this is certainly true for $k\le0$, where we have $x_k\equiv0$.
Thus, with $\eps<\frac12$ for $k\ge1$ and induction on $k$ we get
\begin{eqnarray*}\textstyle
L_k &=& 2|x_{k-1}|+7|x_{k-2}|+7|x_{k-3}|+2|x_{k-4}|\\
    &<&\textstyle
        (\frac6\eps)^{k-1}
[\frac23+\frac73\frac\eps6+
         \frac73(\frac\eps6)^2+
         \frac23(\frac\eps6)^3]\ \ <\ \ 
                 (\frac6\eps)^{k-1}\ =\ b_k
\end{eqnarray*}
and thus the right-hand side in~(\ref{eq:kth-ineq}) is 
always strictly positive, and we also get the inequality 
$|x_k|<\frac1\eps(b_k+L_k)=\frac2\eps(\frac6\eps)^{k-1}<\frac13(\frac6\eps)^k$.
\end{proof}

\subsection{Strictly preserved faces}
\label{subsec:strictlypreservedfaces}

In the following, we are considering an arbitrary 
$m$-dimensional polytope $P\subset\R^m$,
but of course you should think of~$P=D_m^\eps$,
the polytope that we will want to apply this to.

The nontrivial faces $G\subseteq P$ of such a polytope are
defined by linear functions: A nonzero linear function $x\mapsto c^tx$
\emph{defines} the face $G\subseteq P$ if $G$ consists
of the points of~$P$ for which the value $c^tx$
is maximal, that is, if 
\[
G\ =\ \{x\in P:c^tx=c_0\} \ =\ P\cap H,
\]
where $c_0=\max\{c^tx:x\in P\}$, and 
where $H=\{x\in\R^m:c^tx=x_0\}$ is a hyperplane.

\begin{figure}[ht]
\begin{center}
\input{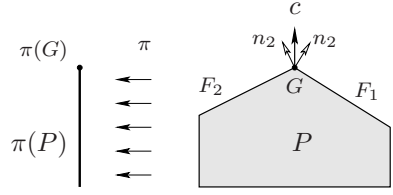}
\end{center}
\vspace{-2mm}
\caption{The normal vector $c$ for the face~$G$ can be written
as a positive combination of the normal vectors of the faces
$F_1,F_2$ that contain~$G$.
The face $G\subset P$ is strictly preserved by the projection
of~$P$ to the second coordinate.}
\label{figure:faceGpreserved}
\end{figure}

Given a face $G$, how do we find a linear functional $c^tx$ that
defines it? It is easy to check (see \cite[Lect.~2]{Z35} for proofs,
and Figure~\ref{figure:faceGpreserved} for intuition)
that $c$ defines $G$ if and only if it is a linear combination,
with positive coefficients, of facet normal vectors $n_F^{}$ of those facets
$F\subset P$ that contain~$G$.

In particular, the affine hull of~$G$, $\aff G$, is the intersection
of all the hyperplanes spanned by the facets $F$ that contain~$G$:
\[
\aff G\ =\ \{x\in\R^m: n_F^{}{}^tx=\max\textrm{ for all facets }F\supseteq G\}.
\]

Now we look at a projection of~$P$, that is, 
we look at a surjective linear map $\pi:\R^m\rightarrow\R^d$.
The image $\pi(P)$ is then a $d$-dimensional polytope,
and the faces of~$\pi(P)$ are all induced by faces of~$P$:
If $\bar G\subset\pi(P)$ is a face of~$\pi(P)$, then
$\pi^{-1}(\bar G)$ is a unique face of~$P$.

Indeed, the faces of the projection, $\bar G\subset \pi(P)$,
thus correspond to the faces of~$P$ that are defined by hyperplanes
that are parallel to the kernel of the projection.
Equivalently, the faces of~$\pi(P)$ correspond to those
faces of~$P$ that are defined by normal vectors that are orthogonal
to the projection.

However, in general the face $\pi^{-1}(\bar G)$ is not the only
face that projects to~$G$, and in general it will have a higher
dimension than $G$, and it will have faces that do not 
project to faces of~$\pi(P)$. (See Figure~\ref{figure:faceGpreserved}
for examples.)
Thus, we single out a very specific, nice situation,
where this does not happen: $G$ will map to a face $\pi(G)\subseteq\pi(G)$
of the same dimension as $G$, and all the faces of~$G$ map to
the faces of~$\pi(G)$.

\begin{definition}[Strictly preserved faces]
Let $\pi:P\rightarrow \pi(P)$ be a polytope projection.
A nontrivial face $G\subset P$ is \emph{strictly preserved}
by the projection if $\pi(G)$ is a face of~$\pi(P)$,
with $G=\pi^{-1}(\pi(G))$, and such that the map 
$G\rightarrow \pi(G)$ is injective.
\end{definition}

One can work out linear algebra conditions that characterize
faces $G$ that are strictly preserved by a projection (see \cite{Z97}):
We need that the normal vectors $n_F$ to the facets $F$ that contain $G$,
after projection to the kernel (or to a fiber) of the projection
do span this fiber positively, that is, the projected vertors
have to span the fiber, and they have to be linearly dependent with
positive coefficients.

Here we want to apply this only in a very specific situation,
namely for an orthogonal projection ``to the last $k$ coordinates,''
that is, for a projection $\pi:\R^m\rightarrow\R^k$
given by $x=(x',x'')\mapsto x''$, where
$x''$ denotes the last $k$ coordinates of~$x$, and $x'$
denotes the first $m-k$ coordinates.
For this situation, the characterization of strictly preserved faces
boils down to the following.

\begin{lemma}
  Let $P\subset\R^m$ be an $m$-dimensional polytope, and
  let $\pi:\R^m\rightarrow\R^k$, $(x',x'')\mapsto x''$
  be the projection to the last $k$ coordinates, which
  maps $P$ to the $k$-polytope $\pi(P)$.
 
  Then a nontrivial face $G\subset P$ is strictly preserved
  by the projection if and only if the facet normals $n_F$ 
  to the facets $F\subset P$ that contain~$G$ satisfy the following
  two conditions:
    Their restrictions $n_F'\in\R^{m-k}$ to the first $m-k$ coordinates
  \begin{compactitem}[~$\bullet$~]
  \item
    must be positively dependent, that is, they must satisfy a
    linear relation of the form $\sum\limits_{F\supset G}\lambda_F n_F'=0$
    with real coefficients $\lambda_F>0$, 
  \item
    and they have full rank, that is, the vectors
    $n_F'k$ span $\R^{m-k}$.
  \end{compactitem}
\end{lemma}

\subsection{Positive row dependencies for the matrix \boldmath$A'_m$}
\label{subsec:positiverowdependencies}

Our aim in the following will be to prove that lots 
of faces of the deformed cube $C_m^\eps$ constructed in
Section~\ref{subsec:constructdeformedcube} survive the 
projection $\pi:\R^m\rightarrow\R^4$ to the last four coordinates;
in particular, we want to see that all the $2$-faces of
the surface $Q_m^\eps$ survive the projection.

In view of the criteria just discussed, we have to verify
that the corresponding rows of the matrix from (\ref{eq:defcube}),
after deletion of the last four components, are positively dependent
and spanning.
This may seem a bit tricky because of the $\eps$ coordinates
around, and because we have to treat lots of different faces,
and thus choices of rows. However, it turns out to be 
surprisingly easy.

We start with the matrix $A_m'\in\R^{m\times(m-4)}$,
\[
A_m'\ :=\ 
\left(\begin{array}{r@{\hspace{5pt}}r@{\hspace{5pt}}r@{\hspace{5pt}}r@{\hspace{5pt}}r@{\hspace{5pt}}r}
  0 \\
  2 &  0 \\
 -7 &  2 &  0 \\
  7 & -7 &  2 &  0 \\
 -2 &  7 & -7 &  2 & \cdot\\
    & -2 &  7 & -7 &\cdot& 0\\
    &    & -2 & \cdot  & -7 &  2 \\
    &    &    & \cdot  &  7 & -7 \\
    &    &    &        & -2 &  7 \\
    &    &    &        &    & -2 
  \end{array}\right).
\]
This is the matrix that you get from the left-hand side matrix of
(\ref{eq:defcube})
if you put $\eps$ to zero, and if you delete the last
four coordinates in each row.

The vectors
\begin{eqnarray*}
  (1,0,0,&\ldots&,0)\\
  (1,1,1,&\ldots&,1)\\
  (1,2,4,&\ldots&,2^n)\\ 
  (1,\tfrac12,\tfrac14,&\ldots&,\tfrac1{2^n})
\end{eqnarray*}
lie in the kernel of this matrix, that is, they describe row dependencies.
Indeed, the coefficients $(2,-7,7,-2)$ that appear in the columns
of~$A_m^\eps$, and hence of~$A_m'$, have been chosen exactly to make this
true.

In particular, the rows of~$A_m'$ are positively dependent
with the coefficient~$0$ for the first row, and
\[
(2^{i-t}-1)(1-2^{t+1-i})\ =\ 
2^{-t}2^i + 2^{t+1}\tfrac1{2^i}-3\qquad
\textrm{for}\quad2\le i\le m-4.
\]
These coefficients are positive, except for the coefficients
for $i=0,t,t+1$, which are zero.
Thus, if we delete the first, $t$-th and $(t+1)$-st row from
$A_m'$, the remaining $m-3$ rows are positively dependent.
Moreover, the remaining $m-3$ rows span $\R^{m-3}$, as one
sees by inspection of~$A_m'$: The rows $2,\dots,t-1$ 
have the same span as the first $t-2$ unit vectors $e_1,\dots,e_{t-2}$,
since the corresponding submatrix has lower-triangular form with
diagonal entries $+2$,
and the rows numbered $t+2,\dots,m$ together have the same span
as $e_{t-2},\dots,e_{m-4}$,
due to a corresponding upper-triangular submatrix with
diagonal entries~$-2$.

So the $m-3$ rows from $A_m'$ corresponding to the index set
$[n]\setminus\{1,t,t+1\}$ are positively dependent and spanning,
for $1<t<n$.
In particular, this is true for the rows with index set
$[n]\setminus\{t,t+1\}$ for $1\le t<n$ as well as for the
rows given by $[n]\setminus\{1,n\}$.
That is, if we delete any two cyclically-adjacent rows from $A_m'$,
then the remaining rows are positively dependent and spanning.
Moreover, the property of a vector configuration to
be ``positively dependent and spanning'' is
stable under sufficiently small perturbations:
Thus if we delete the last four columns, the first row,
and any two adjacent rows from~$A_m^\eps$,
then the rows of the resulting matrix will
be positively dependent, and spanning.
Thus we have proved the following result.

\begin{proposition}\label{prop:embedQm}
For sufficiently small $\eps>0$, 
the projection $\pi:\R^m\rightarrow\R^4$ 
yields a polyhedral embedding of the surface $Q_m^\eps$
in~$\R^4$, as part of the boundary complex of the 
polytope $\pi(D_m^\eps)$.
\end{proposition}

\subsection{Completion of the construction, via Schlegel diagrams}

In the last section, we have constructed 
a $4$-dimensional polytope 
\[\bar P_m\ :=\ \pi(D_m^\eps)\subset\R^4\]
as the projection of an $m$-cube.
One can quite easily prove that the projection is 
in sufficiently general position with respect to the $m$-cube,
so the resulting $4$-polytope is \emph{cubical}: All its facets
are combinatorial cubes.

Moreover, all the vertices and edges of this polytope are
induced from the $m$-cube: We have constructed 
\emph{neighborly cubical} $4$-polytopes.
(Indeed, they are very closely related to the 
neighborly cubical $4$-polytopes as first constructed by
Joswig \& Ziegler \cite{Z62}.)

The boundary complex of any $4$-polytope may be visualized
in terms of a Schlegel diagram (see \cite[Lect.~5]{Z35}):
By stereographic projection from a point that is very close
to a facet $F_0\subset\bar P_m$, we obtain a 
polytopal complex $\mathcal D(\bar P_m,F_0)$ that
faithfully represents all the faces of~$\bar P_m$, 
except for $F_0$ and $\bar P_m$ itself.
Hence we have arrived at the goal of our construction.

\begin{theorem}
  For $m\ge3$, there is a polyhedral realization of the
  surface $Q_m$, the ``mirror complex of an $m$-gon,''
  in~$\R^3$.

  For $m\ge4$ such a realization may be found
  as a subcomplex of
\[
\mathcal D(\pi(D_m^\eps),F_0),
\]
  the Schlegel diagram (with respect to an arbitrary
  facet $F_0$) of a projection of the deformed $m$-cube 
  $D_m^\eps\subset\R^m$ (with sufficiently small $\eps$)
  to the last $4$ coordinates.
\end{theorem}

\begin{figure}[ht]
\begin{center}
\input{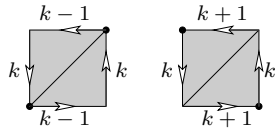}
\end{center}
\caption{The triangulation of~$Q_m$ described above. 
Here we assume that $k$ is even. 
The black vertices are the ones with an even sum of coordinates.}
\label{figure:Q_m-triangulate}
\end{figure}

Thus we have obtained quadrilateral surfaces, polyhedrally
realized in~$\R^3$, of remarkably high genus.
If you prefer to have triangulated surfaces, you may of
course further triangulate the surfaces just obtained,
without introduction of new vertices.
This yields a simplicial surface embedded in~$\R^3$,
with $f$-vector
\[
(2^m, 3m2^{m-2}, m2^{m-1}).
\]
For even $m\ge4$ this may be done in such a
way that the resulting surface has all vertex degrees equal
(to~$\frac32m$): to achieve this, triangulate the
faces with fractional coordinates $k-1$ and $k$ by using the
diagonal between the even-sum vertices if $k$ is even,
and the diagonal between the odd-sum vertices if $k$ is odd.
(Figure~\ref{figure:Q_m-triangulate} indicates how two adjacent
quadrilateral faces are triangulated by this rule.)

In other words, this yields \emph{equivelar} triangulated surfaces of high
genus, which is what McMullen et al.\ were after in~\cite{mcmullen83:_polyh_e}.

Let's finally note that this construction has lots of interesting
components that may be further analyzed, varied, and extended.
Thus a lot remains to be done, and further questions abound.
To note just a few aspects briefly:
\begin{compactitem}[~$\bullet$~]
  \item
Give explicit bounds for some $\eps>0$ that is ``small enough''
for Proposition~\ref{prop:embedQm}.
  \item
Are the neighborly cubical $4$-polytopes 
constructed here combinatorially equivalent to those
obtained by Joswig \& Ziegler in \cite{Z62}?
  \item 
There are higher-dimensional analogues of this: So, extend
the construction as given here in order to get
neighborly cubical $d$-polytopes, with the $(\frac d2-1)$-skeleton
of the $N$-cube, for $N\ge d\ge2$. (Compare \cite{Z62}.)
  \item
Extend this to surfaces that you get as ``mirror complexes''
in products of polygons, rather than just $m$-cubes (which are
products of quadrilaterals, for even~$m$).
\end{compactitem}
See Ziegler \cite{Z97} and Schröder \cite{schroeder04} for work
and ideas related to these questions.
\vfill

\subsection*{Acknowledgements}

Thanks to Michael Joswig, Frank Lutz, Raman Sanyal, Thilo Schröder,
and in particular for everyone at the Oberwolfach Seminar
for interesting and helpful discussions, hints, and comments.
Thanks to Torsten Heldmann for~$\infty$.
\eject

\bibliographystyle{siam}
\bibliography{../../bib/POLYref,../../bib/reference,../../bib/POLYref2,../../bib/ziegler,../../bib/topg}
\end{document}